\documentclass[10pt]{amsart}
\usepackage{epsfig}

\newcommand{\abs}[1]{\left| #1 \right|}
\newcommand{\Prob}[1]{\Pr \left[ #1 \right]}
\newcommand{\mean}[1]{{\mathrm E} \left[#1 \right]}
\newcommand{\variance}[1]{{\mathrm{Var}} \left[#1 \right]}
\newcommand{\sumprob}[2]{\Prob{S_{#1} = #2}}

\newcommand{\ninf}{n \rightarrow \infty}
\newcommand{\littleo}[1]{ o \left( #1 \right)}
\newcommand{\bigO}[1]{ O \left( #1 \right)}
\newcommand{\bigOmega}[1]{ \Omega \left( #1 \right)}
\newcommand{\littleofrac}[2]{\littleo{\frac{#1}{#2}}}
\newcommand{\bigOfrac}[2]{\bigO{\frac{#1}{#2}}}

\newcommand{\bigTheta}[1]{\Theta \left( #1\right)}
\newcommand{\intfloor}[1]{\left\lfloor #1 \right\rfloor}

\newcommand{\thetafun}[2]{\theta_{#1} \left( #2 \right)}
\newcommand{\Z}{\ensuremath{\mathbb Z}}
\newcommand{\ZnZ}{\ensuremath{ \Z/n\Z}}
\newcommand{\R}{\ensuremath{\mathbb R}}
\newcommand{\RZ}{\ensuremath{\R/\Z}}
\newcommand{\tvdist}[2]{\left\| #1- #2 \right\|_{TV}}
\newcommand{\nm}{n+(\mu-1)m}

\newcommand{\cntfrac}{\frac{t-(c_0-1)r}{n+(c_0-1)m}}

\title{A local limit theorem for a family of non-reversible Markov chains}
\author{Elizabeth L. Wilmer}
\date{\today}
\address{Department of
Mathematics\\ Oberlin College \\Oberlin, OH 44074}
\thanks{Research supported by the AT\&T Bell Laboratories Graduate Research
Program for Women and the AT\&T Laboratories Fellowship Program.}
\email{elizabeth.wilmer@oberlin.edu}
\newtheorem{Theorem}{Theorem}
\newtheorem{Proposition}{Proposition}
\newtheorem{Lemma}{Lemma}

\theoremstyle{definition}

\theoremstyle{remark}
\newtheorem*{Remark}{Remark}

\begin{document}
\begin{abstract}
By proving a local limit theorem for higher-order transitions, we determine 
the time required for 
\textit{necklace chains} to be close to stationarity. Because necklace chains,
built by arranging identical smaller chains around a directed cycle,
are not reversible, have little symmetry, do not have uniform stationary
distributions, and can be nearly periodic, prior general bounds on rates of
convergence of Markov chains either do not apply or give poor bounds. Necklace
chains can serve as test cases for future techniques for bounding rates of
convergence.  
\end{abstract}
\maketitle
\section{Introduction}
Determining the rate of convergence to stationarity of finite ergodic Markov
chains is a problem of both theoretical and practical interest, with
applications to sampling and estimation problems. It is also
a difficult problem. Over the last 20 years, many techniques have been
developed  for bounding convergence
behavior; see Aldous and Fill~\cite{aldfill},
Behrends~\cite{behrends00}, or Lov\'{a}sz~\cite{lovasz96}
for surveys.
Coupling, strong stationary time, and finite Fourier analysis arguments
all exploit chain symmetry. Second-largest eigenvalue techniques and
inequalities inspired by differential geometry generally require reversibility;
symmetrized versions of those bounds, due to Diaconis
and Saloff-Coste~\cite{diasal96b} and Fill~\cite{fill91}, can be applied to 
non-reversible chains, but only those with strong types of aperiodicity.

Despite this plethora of methods, the time required to be close to stationary
has been determined precisely---with, say, a correct leading term 
constant---only for certain families of chains (Diaconis~\cite{dia96} gives a
survey focussed on chains whose time to be near stationarity displays sharp
cutoffs).  The current work adds the family of
\textit{necklace chains} to the list.  Necklace chains have little symmetry, 
have non-uniform stationary distributions, are not reversible, and can be
nearly periodic (due to deterministic transitions). Thus existing general
bounds on rate of convergence are difficult to apply to necklaces. We hope our
results will allow necklace chains to serve as test cases for techniques for
bounding rates of convergence, and we provide an example of this utility.

\begin{figure}
\begin{minipage}{1.5in}
(a) \\ \mbox{}\qquad \epsfig{figure=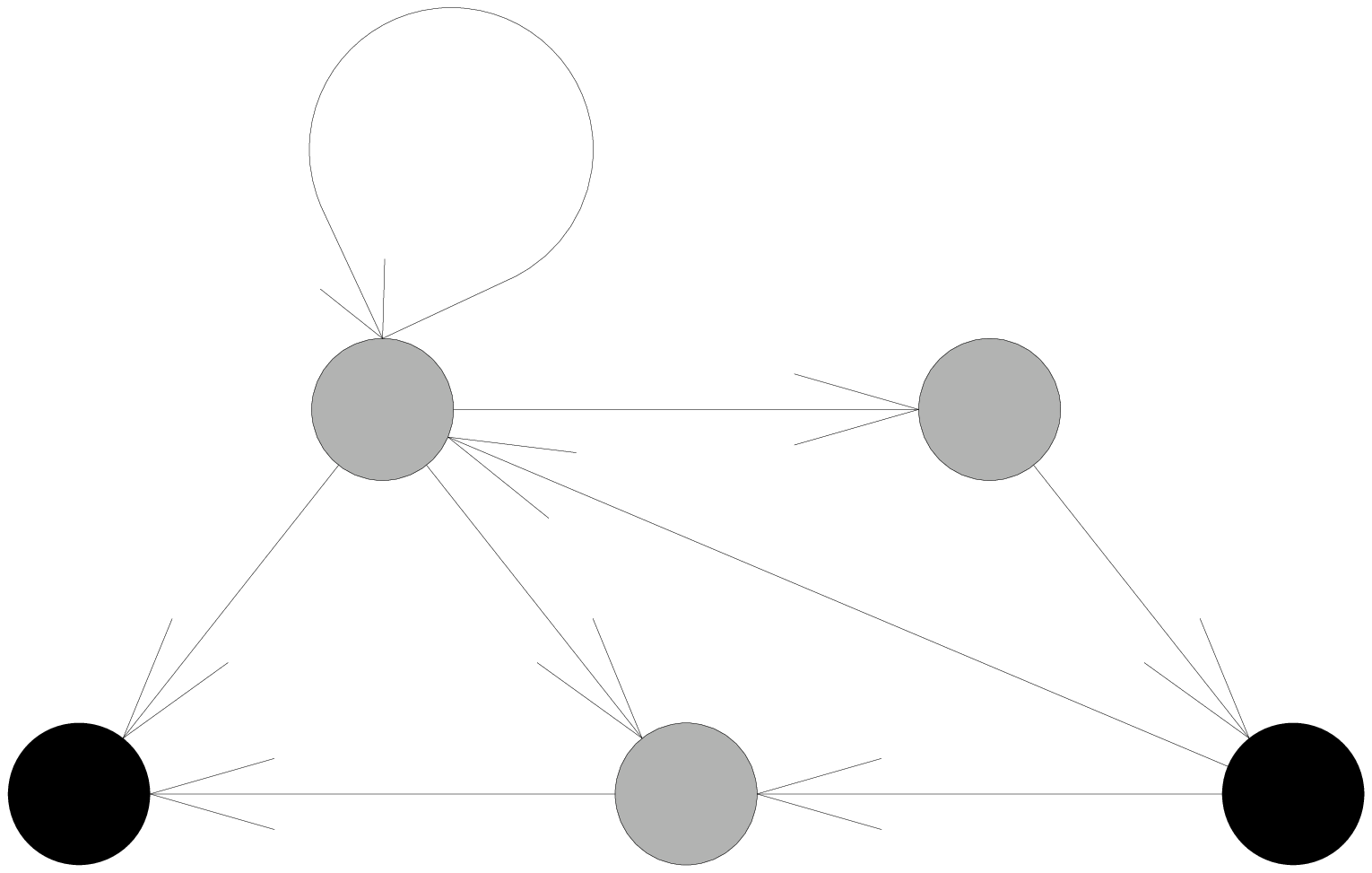,width=30mm}\\[.5in]
(b) \\ \mbox{}\qquad \epsfig{figure=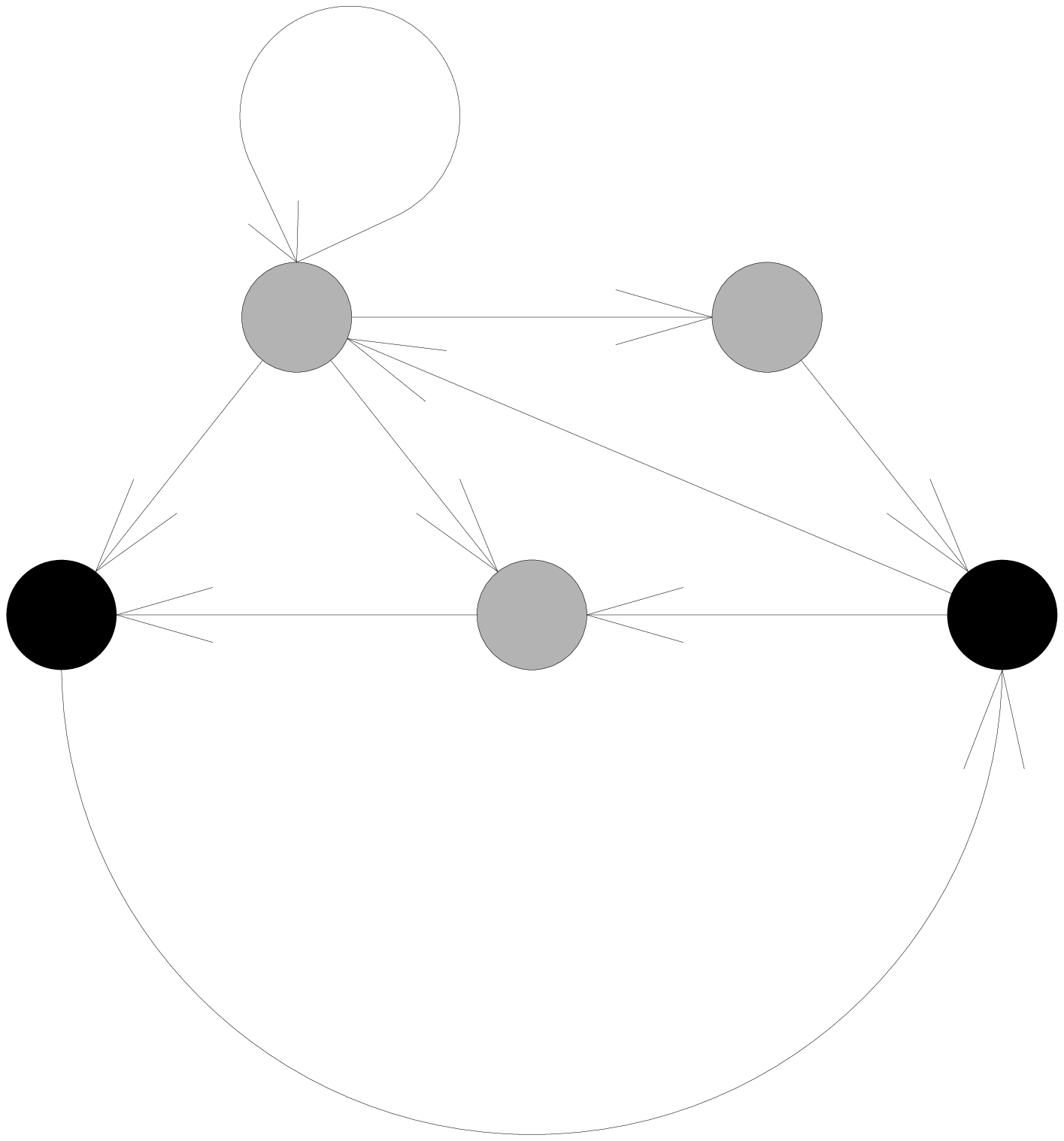,width=30mm}
\end{minipage}
(c)  \begin{minipage}{3in}\mbox{} \qquad \epsfig{figure=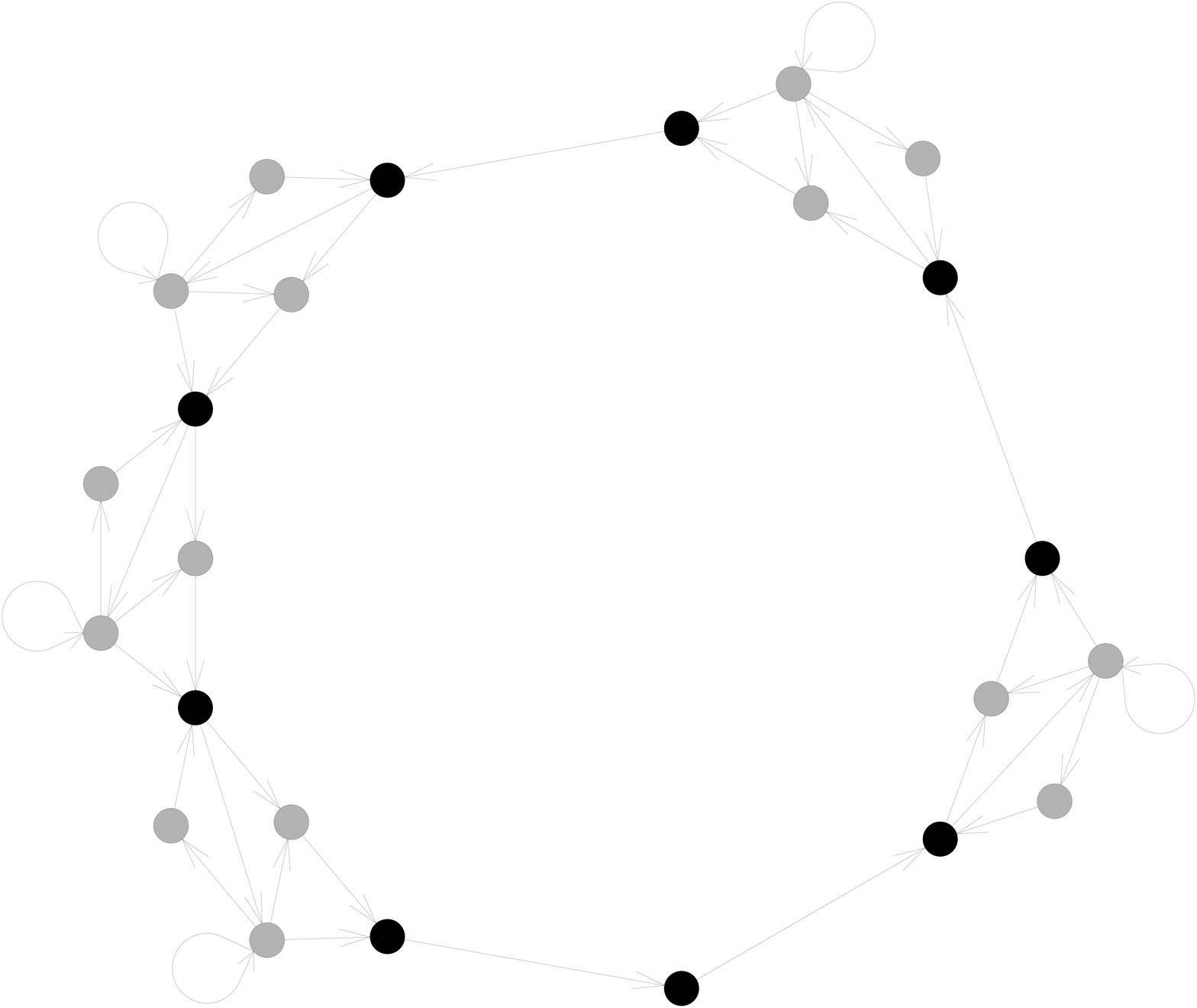,width=70mm}
\end{minipage}
\caption{The underlying graphs of (a)~a bead $B$,
(b)~its closure $\overline{B}$, and (c)~a necklace chain built with
$B$. Here, $n=9$, $m=5$, and $r=(1,0,1,0,1,1,1,0,0)$. Entrance, exit, and link
states are black. 
\label{necklacefig}}
\end{figure}

Necklaces are built from smaller chains.
A \textit{bead} $B$ is a Markov chain with states
$\{0,1,\dots,b\}$ such that
\begin{itemize}
\item The only absorbing state is $b$, and every state lies on some path from
$0$ to $b$. 
\item The set of possible first passage times from $0$ to $b$ 
has minimal span $1$.
\end{itemize}
State $0$ is the {\em entrance state} of $B$
and state $b$ is the {\em exit state}.  The {\em closure}
$\overline{B}$ of the bead $B$ 
has the same  transitions as $B$, except that
$\overline{B}(b,0)=1$.  See Figure \ref{necklacefig}(a,b).

To construct $P_{r}$, a {\em necklace
built with bead $B$ and indicator
vector} $r=(r_0,\dots, r_{n-1})$, start with {\em link states} 
$s_0$, $s_1, \dots, s_{n-1}$.
\begin{itemize}
\item If $r_i=1$,  \textit{there is a bead at position $i$}: attach $s_i$
and $s_{i+1}$ via a bead  isomorphic to $B$
with states $s_{i,0}=s_i$, $s_{i,1}, 
\dots,s_{i,b}=s_{i+1}$.  The \textit{states in
the bead  at position $i$} are $s_{i,0}=s_i$, $s_{i,1}, 
\dots, s_{i,b-1}$. (We exclude
the link state $s_{i,b}=s_{i+1}$ from  the $i$-th bead.)
\item If $r_i=0$,  \textit{there is not a bead at position $i$}: attach $s_i$
and $s_{i+1}$  with a directed  \textit{link edge}, and set
$P_{r}(s_i,s_i+1) =1$.  
\end{itemize}
Indices are taken mod ${n}$, so that $r_{n-1}$ determines
how $s_{n-1}$ is connected to $s_0$. 
See Figure~\ref{necklacefig}(c). 
Let $R_i = \sum_{k=0}^{i-1} r_k$ and let
$ m=R_{n}$ be the total number of
beads.
The Markov chain $P_r$ is always irreducible; it is aperiodic
as long as $m \ge 1$.

In a \textit{family of necklace chains}
$\{P_n\}=\{P_{n}:n\in \mathbb{N}\}$,  $P_n$ has
$n$  link states, has $m=m(n)$ beads ($1 \leq m \leq n$), and
has indicator vector $r=r(n)$. We generally suppress the dependence on
$n$ of these and other chain parameters. 

Our main results are Theorem~\ref{th:llt}, which approximates the
higher--order transitions of  families of necklaces, and Theorem~\ref{th:rate},
which describes the asymptotics of their distance from
stationarity. When rescaled appropriately, families of necklaces behave like
random walk on the cycle $\ZnZ$. 
They converge gradually to stationarity and
have no cutoffs. The time to be near stationarity
depends only on the number of beads, not on their arrangement, and
can range  from $\bigO{|S|^2}$ to $\bigO{|S|^3}$, where $|S|$ is the size of
the state space. 

Let $\phi(x) =\frac{1}{\sqrt{2
\pi}} e^{-{x^2}/{2}} $ be the standard normal density and let 
\[
\thetafun{c}{x} =
\sum_{n \in \Z}\frac{1}{\sqrt{c} } \phi
\left( \frac{n+x}{\sqrt{c}}\right)
\]
be the density at time $c$ of Brownian motion on the circle $\RZ$ of unit
circumference. We  measure distance between distributions in \textit{total
variation}:
\[
\tvdist{\pi(\cdot)}{\pi'(\cdot)} = \sup_A |\pi(A) -\pi'(A)| 
=\frac{1}{2}\sum_s |\pi(s)-\pi'(s)|.
\]
\begin{Theorem}\label{th:llt} 
Let $B$ be a bead with first passage time from entrance to exit of
mean
$\mu$ and variance
$\sigma^2$.  Let $\{P_{n}\}$ be a family of necklace
chains built with bead $B$. Let $\pi_n$ be the stationary 
distribution of $P_n$.
Fix $c>0$  and let 
\[
t =  t(n)= \frac{c(n+(\mu-1)m)^3}{\sigma^2m}+\bigO{n}.
\]
For each $n$, let $s$ be either $s_0$ or a state in a bead at
position
$n-1$  , and let $s'$ be either $s_i$ or a state
in a bead at position $i$. Then as $\ninf$,
\[
P_{n}^{t}(s,s') = \pi_{n}(s')
\thetafun{c}{\frac{t-i-(\mu-1)R_i}{n+(\mu-1)m}}
+\littleofrac{1}{n}.
\]
\end{Theorem}
\begin{Theorem}\label{th:rate} Under the assumptions of
Theorem~\ref{th:llt},  as $\ninf$,
\[
\tvdist{P_n^t(s,\cdot)}{\pi_n(\cdot)} \rightarrow \frac{1}{2}
\int_{\R/\Z}
\left| 
\thetafun{c}{x}-1 \right| dx.
\]
\end{Theorem}

Section~2 collects for reference some standard facts and
computes the stationary distributions and the higher-order transitions of
necklace chains.  In Section~3, we prove Theorems~\ref{th:llt}
and~\ref{th:rate}. The argument applies a  local central limit
theorem due to Petrov, a classical large deviation bound, and Markov
chain identities to the combinatorial expressions for the higher-order
transitions. Section 4 examines several families of necklace chains, including
one (with  well-behaved multiplicative reversibilizations) for which a Nash
inequality of Diaconis and Saloff-Coste~\cite{diasal96b} gives the correct rate
of convergence, up to a constant factor.

\section{Preliminaries}

First we state (for future reference) the specializations of two
standard results on Markov chains to $\overline{B}$. 
See, for example, Aldous and Fill~\cite{aldfill}.
\begin{Proposition} \label{pr:expbound} Let $B$ be a bead, and let $X$ be a
random variable distributed as the first-passage time from $0$ to $b$ in
$B$ (or, equivalently, $\overline{B}$). Then there exist a positive integer
$n_0$ and a constant
$0<\alpha<1$ such that $n>n_0$ implies
$\Pr(X>n)<\alpha^n$.
\end{Proposition}
\begin{Proposition} \label{pr:taboosum}Let $B$ be a bead, let $\pi$ be the
stationary distribution of its closure $\overline{B}$, and let $\mu$ be the
expected first passage time from $0$ to $b$ in $B$. Then, for any state $k
\not=b$ of $B$,
\[
\sum_{a \geq 0} B^a(0,k) =  \pi(k)(\mu+1).
\]
\end{Proposition}  
\begin{proof} This sum is the expected
number of visits to state $k$ when the chain $B$ is started in  $0$ and
run until hitting  $b$, which is equal to the corresponding quantity for
the irreducible chain $\overline{B}$. By a standard result (see,
e.g., Aldous and Fill~\cite[Chapter 2, Lemma 9]{aldfill}), this is
\[
\pi(k) \left(E_0T_b+E_bT_k-E_0T_k \right),
\]
where $E_iT_j$ is the expected passage time between $i$ and $j$ in
$\overline{B}$. But $E_0T_b=\mu$ by definition, while $\overline{B}(b,0)=1$
forces
$E_bT_k-E_0T_k=1$. 
\end{proof}

We now compute the stationary distribution and
higher-order transitions of necklace chains. Proposition
\ref{pr:hot}, which expresses the higher-order
transitions of $P_{n,m}$ in terms of sums of  i.i.d.\ random variables, is
our central combinatorial argument. 

\begin{Proposition}\label{pr:statdist} Let $P_{n,m}$ be a necklace
chain with $n$ link states and $m$ beads of type $B$. Let $\pi$ be
the stationary distribution of $\overline{B}$ and let $\mu=E_0T_b$ in
$\overline{B}$. Then  the stationary distribution
$\pi_{n,m}$ of
$P_{n,m}$ is
\[
\pi_{n,m}(s) = \begin{cases}
\displaystyle{\frac{(\mu+1) \pi(k)}{n+(\mu-1)m}} & \textrm{if } s=s_{i,k}
\textrm{ is in a bead,}\\
\displaystyle{\frac{1}{n+(\mu-1)m}} & \textrm{otherwise.}
\end{cases}
\]  
\end{Proposition}
\begin{proof}
The function 
\[
f(s)=  \begin{cases}
\pi(k) & \textrm{if } s=s_{i,k} \textrm{ is in a
bead,}\\
\pi(b) & \textrm{otherwise}
\end{cases}
\]
satisfies $\sum_{s'} P_{n,m}(s,s') f(s) =f(s')$. Why? Non-backbone states are
only accessible from other states in their own bead, so we are just
verifying stationarity of $\pi$ for $\overline{B}$. Backbone states not in
beads receive and emit a steady flow of $\pi(b)$. Backbone states in beads
receive $\pi(b)$ (which in $\overline{B}$ would come from the
state $b$) from what precedes them, whether bead or backbone state.

To normalize, note that $\sum_s f(s) = (n-m)(1-\pi(b)) + m \pi(b)$. 
Because $\overline{B}(b,0)=1$,  
$E_bT^+_b=\mu+1$ and  $\pi(b)={1}/{(\mu+1)}$.
\end{proof} 

\begin{Proposition}\label{pr:hot} 
Let $P_{n,m}$ be a necklace chain with $n$ link states and $m$
beads of type $B$.  Let $X_1,X_2, \ldots$ be i.i.d.
random variables distributed as the first passage time from $0$ to
$b$ in the bead
$B$, and set $S_j=X_1+X_2+\cdots+X_j$. Then
\begin{align}
\label{eq:linkstate}
P_{n,m}^t(s_0,s) & = 
\begin{cases} \displaystyle{\sum_{j \geq 0}
\sumprob{mj+R_i}{t-i+R_i-(n-m)j}} & \\
\hfill\textrm{when } s=s_i \textrm{ is not in a
bead,}\\
\displaystyle{\sum_{a \geq 0} 
B^a(0,k)\sum_{j \geq 0}  
\sumprob{mj+R_i}{t-a-i+R_i-(n-m)j}} & \\
\hfill\textrm{when } s=s_{i,k} \textrm{ is in a bead, } i\not=0,
\\
\displaystyle{\sum_{a \geq 0} 
B^a(0,k)\sum_{j \geq 0}  
\sumprob{mj+R_i}{t-a-i+R_i-(n-m)j}}+B^t(0,k) & \\
\hfill\textrm{when } s=s_{0,k} \textrm{ is in a bead at position $0$}.
\end{cases}
\intertext{When $s_{n-1,l}$
is in a bead at position $n-1$,} 
\label{eq:arbstate}
P_{n,m}^t(s_{n-1,l},s)& = 
\begin{cases} \displaystyle{\sum_{a \geq 0}
B^a(l,b)P_{n,m}^{t-a}(s_0,s)}  \phantom{mmmmmmmmmmmmmmmmmmm}& \\
\mbox{}\hfill \textrm{when } s \textrm{is not in the
bead at position $n-1$,}\\
\displaystyle{ \sum_{a \geq 0}
B^a(l,b)P_{n,m}^{t-a}(s_0,s)+ B^t(l,k)} &  \\
\hfill\textrm{when } 
s=s_{n-1,k} \textrm{ is in the
bead at position $n-1$}.
\end{cases}
\end{align}
\end{Proposition}
\begin{proof}
In the first case of~(\ref{eq:linkstate}), $j$ counts the number of times the
random walk has gone around the entire chain. By time $t$, the walk has
gone through
$mj+R_i$  beads and $(n-m)j+(i-R_i)$ deterministic steps.
Thus it is   at  $s_i$ exactly when  $S_{mj+R_i}+(n-m)j+(i-R_i)=t$.

For the second case of~(\ref{eq:linkstate}), let $x_t$ be the state occupied at
time $t$, and consider a trajectory of length
$t$ from
$s_0$ to
$s_{i,k}$. 
By the same reasoning as for the first case, 
\begin{multline*}
 \Pr[x_t=s_i \text{ and $x_{t-1}$ is 
not in the bead at position $i$} | x_0=s_0] \\
=\sum_{j \geq 0} \sumprob{mj+R(i)}{t-i+R(i)-(n-m)j}.
\end{multline*}
Let $t-a$ be the time of the last arrival
at $s_i=s_{i,0}$ from the $i-1$-st bead before time~$t$.  In order
to stay in the $i$-th bead until arriving at $s_{i,k}$ at time $t$,
we must avoid the state
$s_{i+1}=s_{i,b}$. 

In the third case, it is possible that the chain never
leaves the bead at position 0. 

For~(\ref{eq:arbstate}), note that when $s$ is not in the bead at position
$n-1$, then we must pass through $s_{n-1,b}=s_0$ on the way to $s$. When
the target state is in the same bead, there are also 
paths that stay within a single bead. 
\end{proof}

\section{Proofs of Limit Theorems}
The key to proving Theorem~\ref{th:llt}, and thus
Theorem~\ref{th:rate}, is Lemma~\ref{le:approx}, an
approximation of the sums appearing in Proposition~\ref{pr:hot}.
While the terms of this sum are simply probabilities of sums of i.i.d.
random variables taking on particular values, each term is taken
from a  different distribution.  We approximate the large terms  with a local
central limit theorem due to Petrov. The largest terms,
which are very close to those taken from a single distribution,
build a theta function. The rest of the terms are negligible---most are handled
by the local central limit theorem, while a Chernoff-type large deviation bound
covers those terms for which only a few random variables are added.

\begin{Lemma}
\label{le:approx}
Let $X_1,X_2, \ldots$ be i.i.d. positive--integer--valued
random variables with maximal span 1 such that 
for some $\alpha<1$ and $n_0$, $n>n_0$ implies
$\Pr[X_i\geq n]\leq\alpha^n$. 
Let 
$\mu= \mean{X_i}$,
$\sigma^2 = \variance{X_i}$, and 
 $S_n=X_1 + \cdots + X_n$. 

Fix $c>0$. Let $m=m(n)$, $r=r(n)$ and
$t=t(n)$ be positive-integer-valued functions such that $1 \leq m \leq n$, 
$ 0 \leq r \leq n$, 
and
\[
t=t(n) = \frac{c(\nm)^3}{\sigma^2 m}+ \bigO{n}.
\]
Then, as $n \rightarrow \infty$ and uniformly in $r$,
\[
\sum_{j \geq 0} \sumprob{mj+r}{t+r-(n-m)j} = 
\frac{1}{\nm} \thetafun{c}{\frac{t-(\mu-1) r}{\nm}}+ \littleofrac{1}{n}.
\]
\end{Lemma}

\begin{proof}[Proof of Lemma~\ref{le:approx}]
First we use the given bound on the tail of the distribution of the $X_i$'s to
build a simple Chernoff-style large deviation bound (see
Chernoff~\cite{chernoff52} or Chapter 2 of Janson, Luczak, and
Ruci\'{n}ski~\cite{jlr}). We can couple
$X_1$ to a shifted geometric random variable $Y$ such that $X_1 \leq Y$
always and
\[
\Pr[Y_i \geq n+n_0] = 
\begin{cases}
\alpha^n &   {\textrm{if\ }} n \geq 0, \\
1 & \textrm{otherwise.}
\end{cases}
\]
But then, for any $a$,
\begin{align*}
\sumprob{n}{an} & \leq\Pr[e^{tS_n}  \geq e^{tan}]
  \leq \left(\inf_{t>0}
e^{-ta} \mean{e^{tX_1}}\right)^n
 \leq \left(\inf_{t>0}
e^{-ta}\mean{e^{tY}}\right)^n\\
&  =
\left(\frac{(1-\alpha)(a-n_0)^{n_0}}{(a-n_0+1)^{n_0-1}}
\left(\frac{\alpha(a-n_0+1)}{a-n_0}
\right)^{a}\right)^n.
\end{align*}
The first fraction inside the exponential above is a rational function of $a$;
the second decreases exponentially for $a$ sufficiently large. Thus, there
must exist a $c_0>\mu$ and a $k>0$ such that, for $a \geq c_0$, 
\begin{equation}\label{eq:chernoff}
\sumprob{n}{an} \leq e^{-kan}. 
\end{equation}
Because $t+r-(n-m)j \geq c_0(mj+r)$ is equivalent to
$j \leq \cntfrac$,
 inequality~(\ref{eq:chernoff}) now implies
\begin{align}
\label{eq:extail}
\sum_{j=0}^{\intfloor{\cntfrac}} &
\sumprob{mj+r}{t+r-(n-m)j}  \leq 
  \sum_{j=0}^{\intfloor{\cntfrac}}
   e^{-k(t-(n-m)j)} \\
& \phantom{mmmmm}\leq \left( \cntfrac
\right)e^{-k\left(t-(n-m)\left(\cntfrac\right)\right)}
\notag\\
& \phantom{mmmmm}
=\bigO{n^2}e^{-k c_0\left(\frac{t+(n-m)r}{n+(c_0-1)m} \right)} 
 = \bigO{n^2}e^{-\bigOmega{n}} =\littleofrac{1}{n}. \notag
\end{align}

The given tail bound implies that the $X_i$'s have
moments  of all orders, so we may now apply a local central limit theorem for
lattice random variables due to Petrov~\cite{petrov62}, taking only one term in
the asymptotic expansion and specializing to the lattice of integers: as $N
\rightarrow
\infty$,
\begin{equation} \label{eq:lclt}
\sup_{a \in \Z} \left(1+\left|\frac{a-\mu N}{\sigma
\sqrt{N}}\right|^3\right)
\left(\sumprob{N}{a}-\frac{1}{\sigma \sqrt{N}}\phi \left(\frac{a-\mu
N}{\sigma
\sqrt{N}}\right)
\right) =
\littleofrac{1}{\sqrt{N}}.
\end{equation}
Because the underlying random variables, $X_i$, are positive integers, all
terms of our sum for which $mj+r>t+r-(n-m)j$, or equivalently $ j>\frac{t}{n}$,
are zero.  Given the estimate of equation~(\ref{eq:extail}), we 
can restrict our attention to those terms for which 
$\cntfrac< j \leq \frac{t}{n}$. On this range,
\begin{equation} \label{eq:mjest}
mj+r = \bigTheta{n^2}.  
\end{equation}
Let 
\[ 
x_0 =\frac{t-(\mu-1) r}{\nm}.
\] 
and
\[
y_j=\frac{t+r-(n-m)j-\mu(mj+r)}{\sigma \sqrt{mj+r}}=\frac{(\nm)(j-x_0)}{\sigma
\sqrt{mj+r}}.
\]
Because $\cntfrac <j \leq \frac{t}{n}$ implies $y_j=\bigTheta{1} (j-x_0),$
\begin{equation} \label{eq:yjsumest}
\sum_{\cntfrac <j \leq \frac{t}{n} }\frac{1}{1+|y_j|^3} <\infty.
\end{equation}
Combining~(\ref{eq:lclt}),~(\ref{eq:mjest}), and~(\ref{eq:yjsumest}) now yields
\begin{multline*}
\sum_{\cntfrac <j \leq \frac{t}{n} } \sumprob{mj+r}{t+r-(n-m)j} \\
= \sum_{\cntfrac <j \leq \frac{t}{n} } \frac{1}{\sigma \sqrt{mj+r}} \phi 
\left( y_j\right)  + \littleofrac{1}{n}.
\end{multline*}
We first consider the terms of the resulting sum for
which $|j-x_0|<n^{1/4}$. On this range, 
\begin{multline}
\label{eq:recipsqrtest}
\frac{1}{\sigma \sqrt{mj+r}} = \frac{1}{\sigma \sqrt{\frac{c
(\nm)^2}{\sigma^2}+\frac{nr}{\nm} + m(j-x_0)}}\\
 = \frac{1}{\sqrt{c} (\nm) \sqrt{1+\bigOfrac{1}{n^{3/4}} }}
 = \frac{1}{\sqrt{c} (\nm)} +\bigOfrac{1}{n^{7/4}}
\end{multline}
and
\begin{align}
\label{eq:yjest}
y_j & =  \frac{(\nm)(j-x_0)}{\sigma
\sqrt{mj+r}}= 
\frac{j-x_0}{\sqrt{c+\bigOfrac{1}{n^{3/4}} }}  = \frac{j-x_0}{\sqrt{c}}
+\bigOfrac{1}{\sqrt{n}}.
\end{align}
By~(\ref{eq:recipsqrtest}),~(\ref{eq:yjest}) and the boundedness of
the derivative of $\phi$, 
\begin{align*}
\sum_{|j-x_0|<n^{1/4}} \frac{1}{\sigma\sqrt{mj+r}}\phi(y_j) 
& = \sum_{|j-x_0|<n^{1/4}} \left(\frac{1}{\sqrt{c} (\nm)} +\bigOfrac{1}{n^{7/4}}
\right) \\ 
& \phantom{= \sum_{|j-x_0|<n^{1/4}} \sum_{|j-x_0|<n^{1/4}}}
\times \left(\phi \left(
\frac{j-x_0}{\sqrt{c}}\right)+\bigOfrac{1}{\sqrt{n}}\right)\\
& =\frac{1}{\sqrt{c}(\nm)}\sum_{|j-x_0|<n^{1/4}} 
\phi \left(\frac{j-x_0}{\sqrt{c}}\right)+
\bigOfrac{1}{n^{5/4}} \\
& =\frac{\thetafun{c}{x_0}}{\nm}+\littleofrac{1}{n}.
\end{align*}

The remaining terms are those for which $\cntfrac < j
<x_0-n^{1/4}$ or $x_0+n^{1/4} < j < \frac{t}{n}$.
Since there are only $\bigO{n^2}$ such terms, each of which is
$\bigO{\phi(n^{1/4})}$, their sum is certainly $\littleofrac{1}{n}$.
\end{proof}

To complete the proof of Theorem~\ref{th:llt} 
we need only check that
transition probability factors weighting the inner sums in
Proposition~\ref{pr:hot} give the correct the stationary distribution factors
in the final approximation. 

\begin{proof}[Proof of Theorem~\ref{th:llt}] First, note that 
Proposition~\ref{pr:expbound} legitimates applying 
Lemma~\ref{le:approx} to sums arising from 
Proposition~\ref{pr:hot}. 

When $s=s_0$ and $s'=s_i$, where $s_i$ is not in
a bead, applying Lemma~\ref{le:approx} to the sum given by the first case of
Proposition~\ref{pr:hot}(\ref{eq:linkstate}) and recalling
Proposition~\ref{pr:statdist} immediately yields the desired
\begin{align*}
P_{n,m}^t(s_0,s_i)
& = 
\pi_n(s_i)\thetafun{c}{\frac{t-
i-(\mu-1)R_i}{n+(\mu-1)m}}+\littleofrac{1}{n}.
\end{align*}

Next, we consider $s'=s_{i,k}$, a state in a bead at position $i$. Here, the
second case of Proposition~\ref{pr:hot}(\ref{eq:linkstate}) applies. First we
apply Lemma~\ref{le:approx} to each term; by Proposition~\ref{pr:taboosum}, we
can collect the error terms. Then we use the upper bound on
$B^a(0,k)$ implied by Proposition~\ref{pr:expbound} and the boundedness of both
$\theta_c$ and its derivative to truncate the sum, approximate the theta
function factor in the remaining terms by its value when $a=0$, and
de-truncate, all with small enough error. Finally,
Proposition~\ref{pr:taboosum} and Proposition~\ref{pr:statdist}  evaluate the
remaining sum.
\begin{align*}
P_{n,m}^t(s_0,s_{i,k}) 
& = 
\sum_{a \geq 0} \frac{B^a(0,k)}{\nm} 
\thetafun{c}{\frac{t-a-
i-(\mu-1)R_i}{n+(\mu-1)m}}+\littleofrac{1}{n}\\
& = 
\sum_{a=0}^{\sqrt{n}} \frac{B^a(0,k)}{\nm}\left( \thetafun{c}{\frac{t-
i-(\mu-1)R_i}{n+(\mu-1)m}}+\bigOfrac{1}{\sqrt{n}}\right)+\littleofrac{1}{n}\\
& = \frac{\sum_{a \geq 0} B^a(0,k)}{\nm} \thetafun{c}{\frac{t-
i-(\mu-1)R_i}{n+(\mu-1)m}}
 +\littleofrac{1}{n}\\
&=
\pi_n(s_{i,k}) \thetafun{c}{\frac{t-
i-(\mu-1)R_i}{n+(\mu-1)m}}
 +\littleofrac{1}{n}.
\end{align*}
Furthermore, Proposition~\ref{pr:expbound} ensures that for
large $n$ the extra term in the third case of
Proposition~\ref{pr:hot}(\ref{eq:linkstate}) can be absorbed into the error.

We must still consider $s=s_{n-1,l}$ (where $l \not= b$), a state in a bead at
position
$n-1$. Now  Proposition~\ref{pr:hot}(\ref{eq:arbstate}) applies. This time,
the weights sum as $\sum_{a
\geq 0} B^a(l,b)=1$, because the chain $B$ started at $l$ will eventually hit
the absorbing state $b$. Thus the error terms
can be collected. As above, Proposition~\ref{pr:expbound} covers the extra term
in the second case of Proposition~\ref{pr:hot}(\ref{eq:arbstate}). From this
point, we truncate, approximate, de-truncate, and sum as above, obtaining 
\begin{align*}
P_{n,m}^t(s_{n-1,l}, s')& = \sum_{a \geq 0} B^a(l,b)  \pi_n(s')
\thetafun{c}{\frac{t-a-i-(\mu-1)R_i}{n+(\mu-1)m}} + \littleofrac{1}{n}\\
&=
\pi_n(s') \sum_{a=0}^{\sqrt{n}} B^a(l,b) \left( \thetafun{c}{\frac{t-
i-(\mu-1)R_i}{n+(\mu-1)m}}+\bigOfrac{1}{\sqrt{n}}\right)+\littleofrac{1}{n}\\
& = 
\pi_n(s') \sum_{a \geq 0} B^a(l,b) \thetafun{c}{\frac{t-
i-(\mu-1)R_i}{n+(\mu-1)m}}+\littleofrac{1}{n}\\
&=
\pi_n(s') \thetafun{c}{\frac{t-
i-(\mu-1)R_i}{n+(\mu-1)m}}+\littleofrac{1}{n}.
\end{align*}
\end{proof}
We now substitute Theorem~\ref{th:llt} into the $L^1$
expression for the total variation distance and collect terms to obtain a
Riemann sum.
\begin{proof}[Proof of Theorem~\ref{th:rate}] Let $f(x)=|\thetafun{c}{x}-1|$ and
define $z_0, z_1, \dots, z_{n-1} \in \RZ$ by
\[
z_i = \begin{cases}
\displaystyle{\frac{-t}{\nm}} &  i=0\\[3mm]
\displaystyle{z_{i-1}+\frac{1 + (\mu-1)r_i}{\nm}}& 1 \leq i \leq n-1.
\end{cases}
\]
The $z_i$'s cover the entire circle $\RZ$ with spacing of $\bigOfrac{1}{n}$;
the size of the interval between $z_{n-1}$ and $z_0$ is determined by the
presence or absence of a bead at position $n-1$. 

Without loss of
generality, we may assume that
$s=s_0$ or
$s$ is in a bead at position $n-1$.  Then Theorem~\ref{th:llt} implies
\begin{align*}
2 \tvdist{P_n^t(s,\cdot)}{\pi_n(\cdot)} 
& = 
\sum_{s'} \left| P_n^t(s,s') -\pi_n(s')\right|\\
& = 
\sum_{s'} \pi_n(s')
\left|\thetafun{c}{\frac{t-i-(\mu-1)R_i}{\nm}}+\littleo{1}-1
\right| \\
& = 
\sum_{s'} \pi_n(s') f \left(\frac{t-i-(\mu-1)R_i}{\nm} \right) 
+\littleo{1}.\\
\end{align*}
If $s'=s_i$ is not in a bead, then $\pi_n(s')=\frac{1}{\nm}=z_i-z_{i-1}.$ If
there is a bead at position $i$, then Proposition~\ref{pr:statdist}, $r_i=1$,
and $\pi(b)=\frac{1}{\mu+1}$ imply that
\[
\sum_{k=0}^{b-1} \pi_n(s_{i,k}) = \frac{(1-\pi(b))(\mu+1)}{\nm} 
= \frac{\mu}{\nm} = z_i-z_{i-1}.
\]
Now group states by position and recall that $f$ is an even function:
\begin{align*}
2 \tvdist{P_n^t(s,\cdot)}{\pi_n(\cdot)} 
& =
\sum_{i=0}^{n-1} (z_i-z_{i-1}) f(z_i) +\littleo{1} = 
\int_{\R/\Z} f(x) dx +\littleo{1}.
\end{align*}
\end{proof}
\begin{Remark} Let $y^n_0, y^n_1, y^n_2, \dots$ be  trajectory of $P_n$, and
define, for $s \in S$, 
\[
z_t(s) = \frac{i+(\mu-1)R(i)-t}{\nm} \textrm{ whenever $s=s_i$ or $s$ is in a
bead at position
$i$.}
\]
Then 
\[
Y^n(\tau) = z_{\intfloor{\frac{(\nm)^3\tau}{\sigma^2
m}}}\left(y^n_{\intfloor{\frac{(\nm)^3\tau}{\sigma^2 m}}}\right)
\]
is a cadlag random function $[0,1] \rightarrow \R$. Using
a random-change-of-time  argument  parallel to Billingsley's proof of
a functional limit theorem for the renewal process~\cite[pp. 148--50]{bill68},
it can be shown that the sequence
$Y^1, Y^2,
\dots$ converges weakly to Brownian motion on the circle $\RZ$. 
\end{Remark}
\section{Examples}

\subsection*{Fixed number or fixed fraction of beads.}
When $m(n)=m$ is constant, $\{P_{n}\}$ converges to stationarity on a cubic
time scale,
\[
\frac{(n+(\mu-1)m)^3}{\sigma^2 m} = \frac{1}{\sigma^2 m}\left(n^3\right)
+\bigO{n^2}.
\] 
Now let $m(n)=\intfloor{kn}$, where $0 < k \leq 1$. Then $\{P_{n}\}$ converges
to stationarity on the quadratic time scale 
\[
\frac{(n+(\mu-1)\intfloor{kn})^3}{\sigma^2 \intfloor{kn}}
=\frac{(k \mu-k +1)^3}{\sigma^2k}\left(n^2\right) +\bigO{n}. 
\]

\subsection*{A simple bead I: rearranging beads}
\begin{figure}
(a)\qquad \qquad \quad
\raisebox{-.5in}{\epsfig{figure=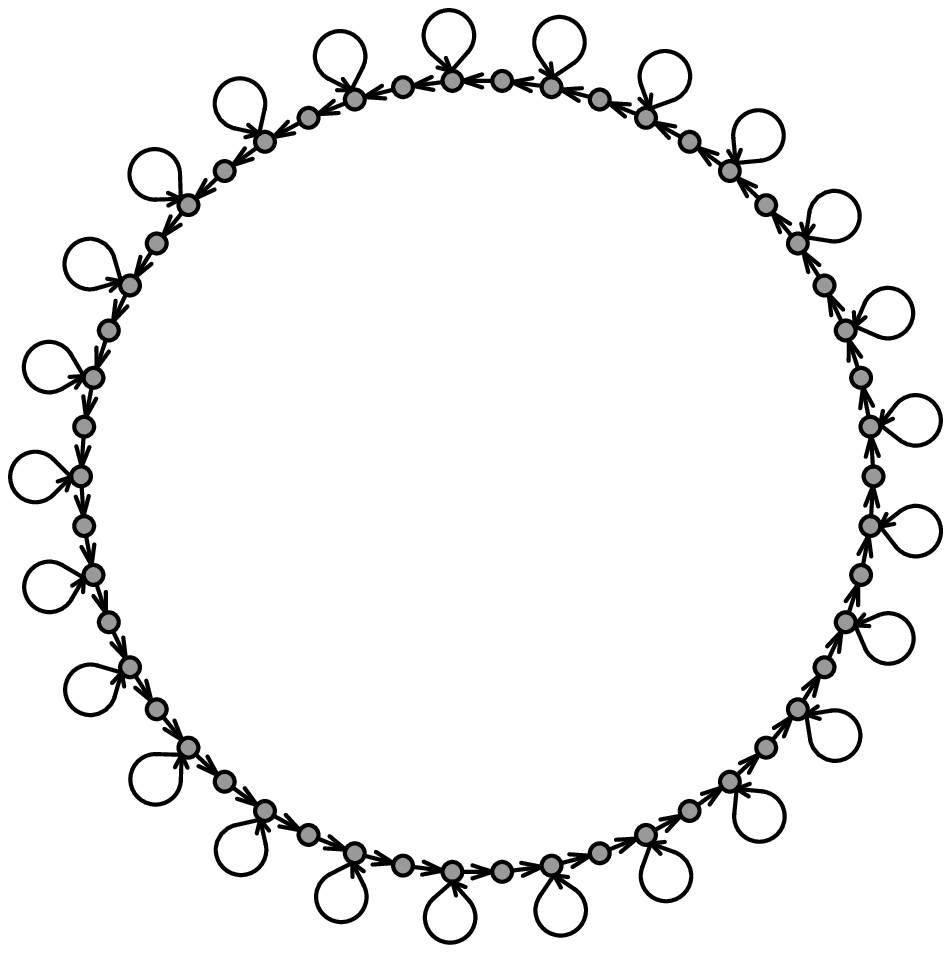,width=45mm}}
\qquad
\raisebox{-.5in}{\epsfig{figure=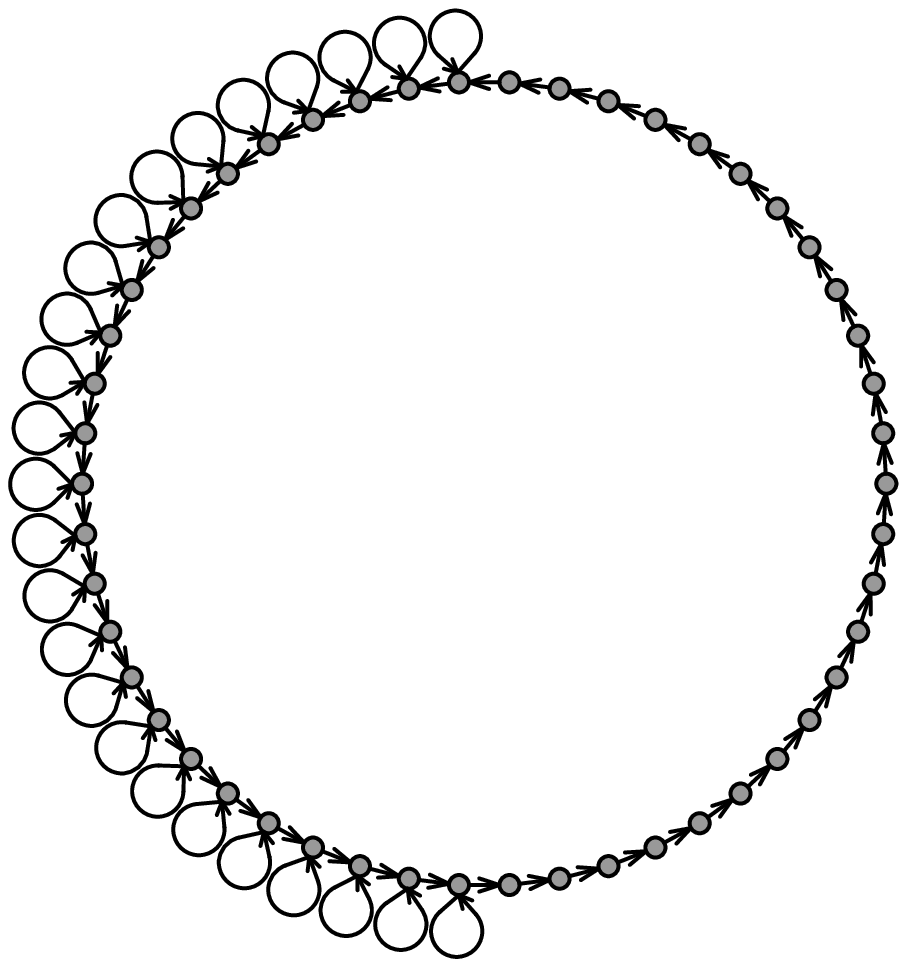,width=45mm}} \qquad \qquad \quad
\\[3mm]
(b) \qquad\raisebox{-.5in}{\epsfig{figure=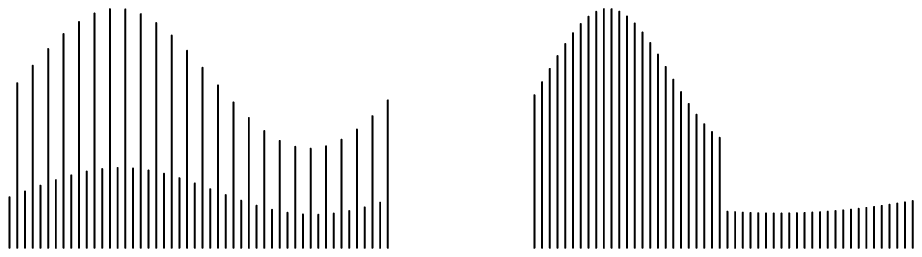,width=110mm}}
\\
(c) \qquad\raisebox{-.5in}{\epsfig{figure=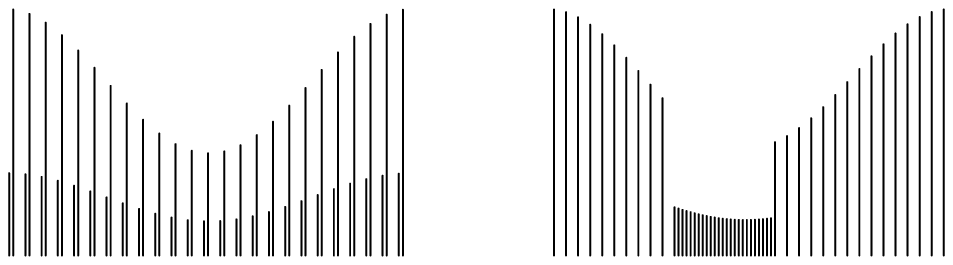,width=110mm}}
\\
(d) \qquad \raisebox{-.5in}{\epsfig{figure=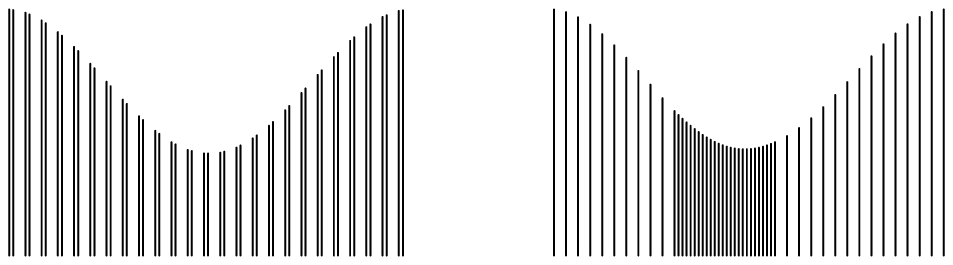,width=110mm}}
\caption{(a)~Underlying graphs of two homogeneous stop/go chains;
 $n=50$, $m=25$, 
$R(i) = \intfloor{(i+1)/2}$ and 
$R(i) = \min(i,n/2).$ (b)~$P_n^t(0,\cdot)$ for the two chains
above;  $p=2/3$ and $t=530$. (c)~$P_n^t(0,\cdot)$
rearranged. (d)~$P_n^t(0,\cdot)$ rearranged and normalized.
\label{illusfig}}
\mbox{} \\
\end{figure}

Consider the bead $B$ with state space $\{0,1\}$
 and transition matrix $\left(\begin{smallmatrix}
p & q \\
0 & 0 
\end{smallmatrix} \right)$.  Then
$B^t(0,1) = qp^{t-1}$, so
$\mu=\frac{1}{q}$ and 
$\sigma^2 = \frac{p}{q^2}$. 
The closure $\overline{B}$ has stationary distribution
$\pi(0)=\frac{1}{1+q}$, $\pi(1)=\frac{q}{1+q}$. Necklace chains built with $B$
are just directed cycles, some of whose states have a hold probability of $p$.
When $\{P_{n}\}$ is a family of such chains, Proposition~\ref{pr:statdist}
implies  
\begin{align*}
\pi_{n}(s) & = \begin{cases}
\frac{1}{qn+pm(n)} &\mbox{if $s$ is in a bead,}\\
\frac{q}{qn+pm(n)}   &\mbox{otherwise,}
 \end{cases} 
\end{align*}
while Theorem~\ref{th:rate} implies  
$\{P_{n}\}$ converges to stationarity on a time scale of
\begin{align}
\label{eq:simplerate}
\frac{(n+(\mu-1)m(n))^3}{\sigma^2m(n)}
& =\frac{(qn+pm(n))^3}{pqm(n)}.
\end{align}
Figure~\ref{illusfig}(a) shows the underlying graphs of two chains built with
this bead: one alternates hold states with forced transitions, the other
groups its holds together.
Since they have the same number of beads, these
chains should evolve at the same rate.
Figure~\ref{illusfig}(b) shows the higher--order transition 
probabilities of these two chains at identical times (about $.08$ times the
time scale of convergence). Notice the  effects of the two-valued
stationary distributions.
Figure~\ref{illusfig}(c) shows the same weights respaced according to 
the points at which Theorem~\ref{th:llt} evaluated $\theta_c$, so that the
value of $P_{n}^t(0,s)$ is attached to
$\frac{t-i-(\mu-1)R_i}{\nm}$.  The two apparent ``curves''
correspond to the two values of the stationary
distribution. 
Figure~\ref{illusfig}(d) normalizes the
weight at $\frac{t-i-(\mu-1)R_i}{\nm}$ by $\pi_{n}(s_i)$. The convergence
claimed in Theorem \ref{th:llt} is now clear.

\subsection*{A simple bead II: optimizing convergence} First, let
$m(n)=n$, so that $\{P_{n}\}$  is a family of walks on $\ZnZ$.
Many techniques for analyzing convergence apply, including
Theorem~\ref{th:rate}: the time scale of convergence is $n^2/pq$, which is
fastest when
$p=q={1}/{2}$. 

Now consider a family  $\{P_{n}\}$
with holds at $\intfloor{kn}$ states, $0 < k < 1$.
One might expect that $p={1}/{2}$ would optimize
convergence again, since
the variance per coin flip is maximized. 
However, this family converges
at time scale 
\[
\frac{(qn+p\intfloor{kn})^3}{pq\intfloor{kn}} =
\frac{(q+pk)^3}{pqk}n^2+\bigO{n}.
\] 
The coefficient of $n^2$ is minimized when
\[
p = \frac{-k+ \sqrt{k^2-k+1}}{1-k}.
\] 
This
optimizing probability is a decreasing function of $k$
which approaches $1$ as $k \rightarrow 0$ and $1/2$ as
$k \rightarrow 1$.

\subsection*{Comparison with other methods}
While it is generally difficult to apply the usual techniques for bounding
rates of convergence to our nearly--periodic examples, there is a
linked pair of examples which are susceptible. For these, we find that a Nash
inequality of Diaconis and Saloff-Coste~\cite{diasal96b} is correct, up to a
constant factor, while the symmetrized second-largest-eigenvalue bound of
Fill~\cite{fill91} is off by a log factor. 

\begin{figure}
\epsfig{figure=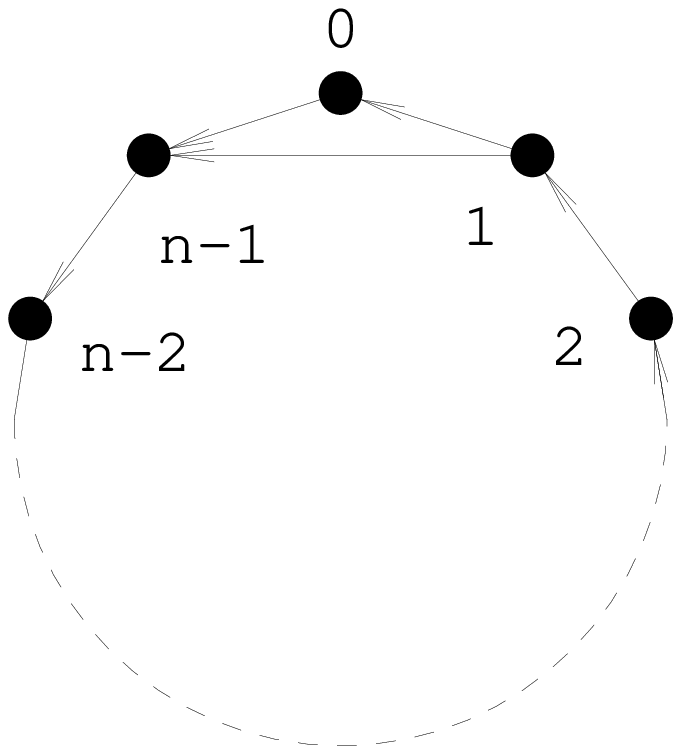,width=35mm} \qquad
\epsfig{figure=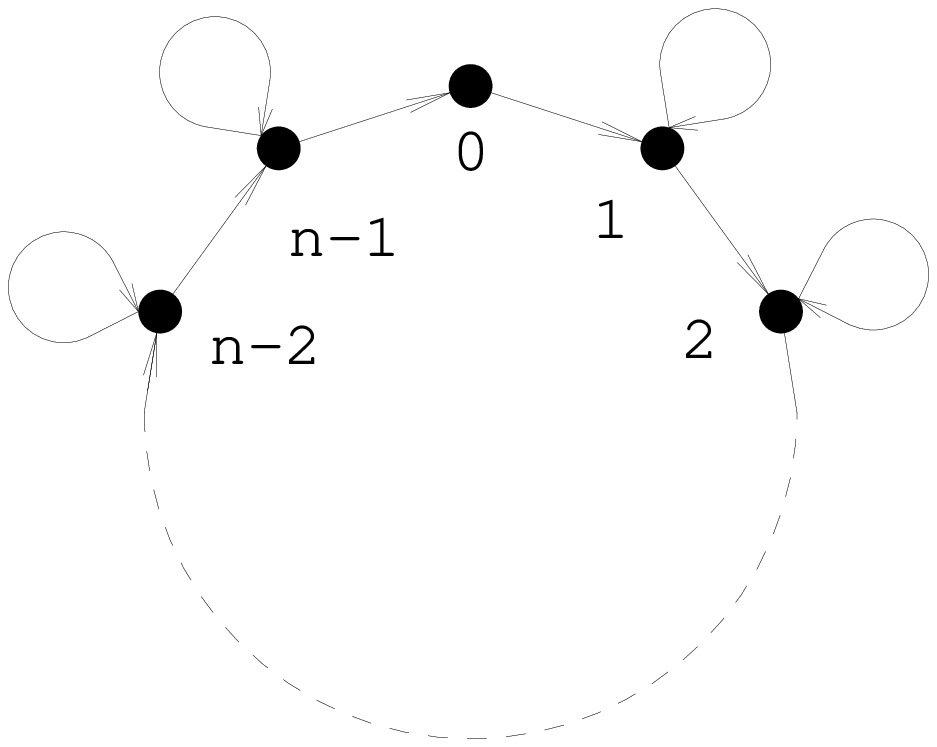,width=35mm}
\caption{The underlying graphs of $P_n$ and $P_n^{n-1}.$ \label{symm1fig}}
\vspace{5mm}
\epsfig{figure=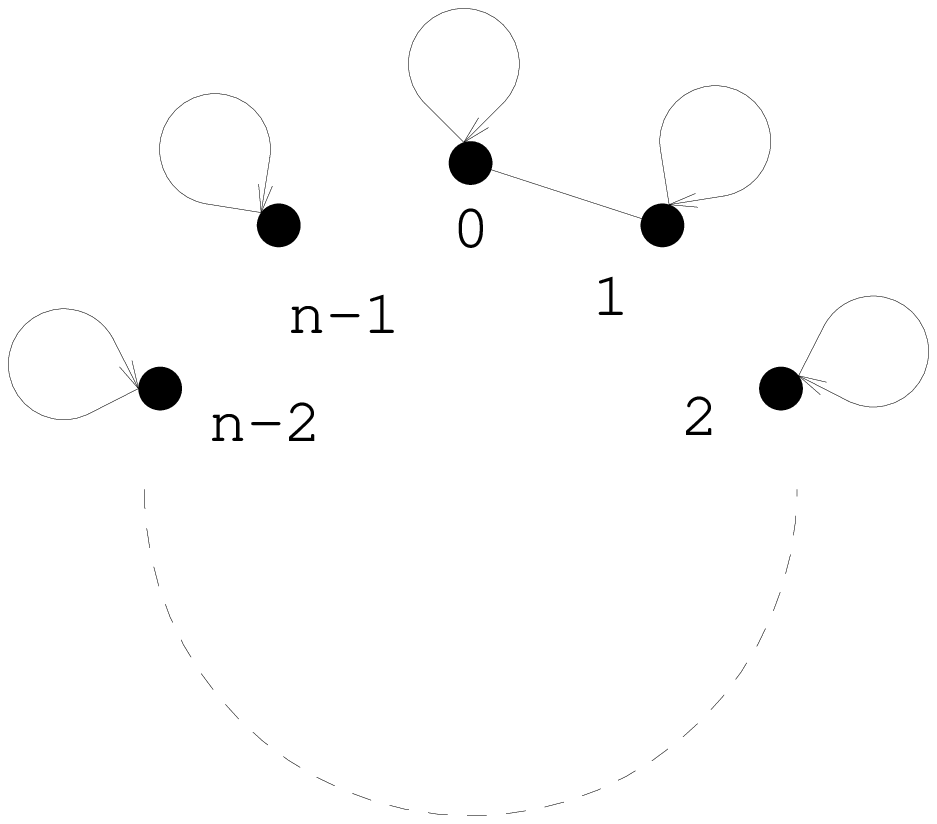,width=35mm} \qquad
\epsfig{figure=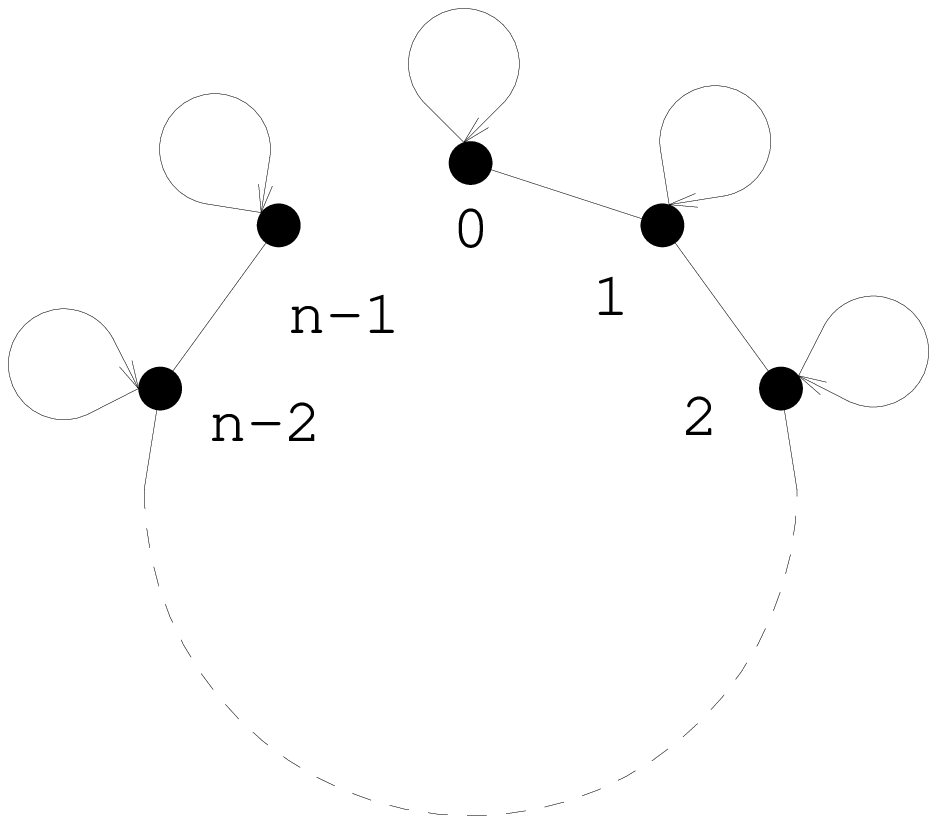,width=35mm}
\caption{The underlying graphs of  $M(P_n)$
 and $K_n=M(P_n^{n-1})$. 
\label{symm2fig}}
\end{figure}

Let
$0 <p<1$,
$q=1-p$, and let 
\[
P_n =
\left[
\begin{array}{ccccc}
0 & 0 & 0 & \dots & 1 \\
 q & 0 & 0 & \dots & p \\ 
0 & 1 & 0 &  \dots & 0 \\ 
\vdots &\ddots  & \ddots & \ddots & \vdots\\ 
0 &  \dots& 0 & 1& 0\\
\end{array} 
\right], \quad
\textrm{so that} \quad
P_n^{n-1} = 
\left[
\begin{array}{cccccc}
 0 & 1 & 0 &\dots & 0 \\
 0 & p & q &  \dots & 0\\ 
0 & 0 & p &  \ddots & \vdots \\
 \vdots & \vdots & \ddots &  \ddots & q \\
 q &0&    \dots & 0 & p 
\end{array}
\right].
\] 
See Figure~\ref{symm1fig}. The chain $P_n$ has
one bead and $n-1$ link states; its bead has $\mu=1+q$ and
$\sigma^2=pq$.  Theorem~\ref{th:rate} gives a convergence time scale of
$(n-p)^3/pq$ for~$\{P_n\}$ .

However, $P_n^{n-1}$ is  also a 
necklace chain, with $n$ link states and $n-1$ simple beads; the
outer cycle has changed direction. By equation~(\ref{eq:simplerate}), the family
$\{ P_n^{n-1}\}$ has convergence time scale 
${(n-p)^3}/{pq(n-1)}$ (consistent
with the time scale for $\{P_n\}$). 

\subsubsection*{Symmetrize and compare}
When $P$ is a Markov chain, let $\overleftarrow{P}$ be the time
reversal of $P$ and $M(P)=P \overleftarrow{P}$ be the
multiplicative symmetrization of $P$. Fill \cite{fill91}
showed that 
\begin{equation}
\label{eq:fill}
\tvdist{P^t(x_0,\cdot)}{\pi(\cdot)} \leq 
\frac{1}{2 \sqrt{ \pi(x_0)}}\beta_1(M(P))^{t/2},
\end{equation}
where $\beta_1(M(P))$ is the second largest eigenvalue of
$M(P)$.

Unfortunately, 
the underlying graph of $M(P_n)$
is almost completely disconnected and thus
$\beta_1(M(P_n))=1$ (see~\cite{diasal96b} for several
similar examples).
However, 
\[
K_n = M(P_n^{n-1}) =
\left[
\begin{array}{ccccc}
    q  & p       & 0       & \dots    & 0 \\
    pq & p^2+q^2 & pq      &  \dots   & 0\\ 
\vdots & \ddots  &  \ddots & \ddots   & \vdots \\
     0 &  \dots  & pq      &  p^2+q^2 & pq \\
     0 & \dots   &  0      & pq       & 1-pq 
\end{array}
\right]
\] 
has underlying graph an
$n$-path (see Figure \ref{symm2fig}). Because 
$K_n$ has non-trivial edge weights, we cannot compute $\beta_1(K_n)$ directly.
Specializing  a comparison result of Diaconis 
and Saloff-Coste~\cite{diasal93b} yields 
\begin{equation}
\label{eq:d-sc}
\beta_1(K_n) \leq 1 - \left(\min_{x}
\frac{\tilde{\pi}_n(x)}{\pi_n(x)}
\right)
\left( \min_{\begin{smallmatrix}K_n(x,y)>0 \\ x \not= y 
\end{smallmatrix}}\,\,
\frac{\pi_n(x)
K_n(x,y)}{\tilde{\pi}_n(x)\tilde{K}_n(x,y)}\right) 
\left( 1 - {\beta_1(\tilde{K}_n})\right),
\end{equation}
whenever $\tilde{K}_n$ is a chain on the same state space as $K_n$ with the
property that $K_n(x,y)>0$ implies $\tilde{K}_n(x,y)>0$; here $\tilde{\pi}_n$ is
the stationary distribution of $\tilde{K}_n$. 

\begin{figure} \hspace{79mm}
$\frac{1-pq}{n-p}$ \hspace{0mm}
$\frac{p^2+q^2}{n-p}$ \hspace{1mm}
\dots
\hspace{1mm}
$\frac{p^2+q^2}{n-p}$ \hspace{0mm}
$\frac{q^2}{n-p}$ \hspace{17mm}\\
 $\pi(x)$ 
\epsfig{figure=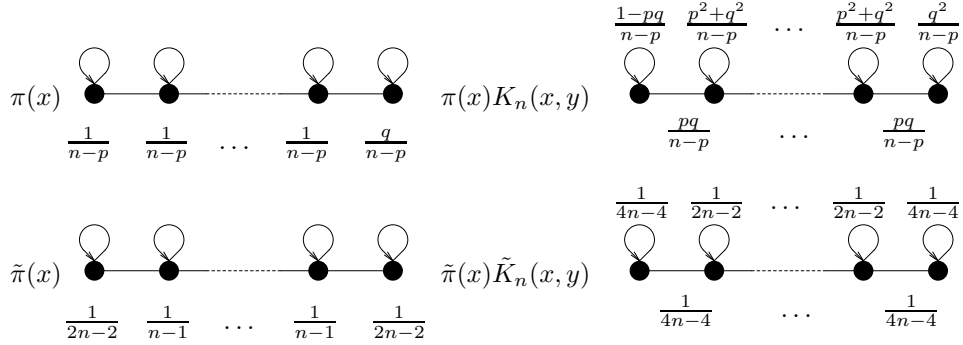,width=46mm}
\hfill
$\pi(x) K_n(x,y)$ \hspace{1mm}
\epsfig{figure=figure5.eps,width=46mm}
\\[.5mm] 
\hspace{1.5mm}
$\frac{1}{n-p}$ 
\hspace{3mm}$ \frac{1}{n-p}$
\hspace{1.5mm} $\dots$ \hspace{1.5mm} $ \frac{1}{n-p}$ 
\hspace{2mm}
$\frac{q}{n-p}$
\hspace{30.5mm} \raisebox{1.5mm}{
$\frac{pq}{n-p}$} \hspace{6.5mm} \raisebox{1.5mm}{$\dots$}
\hspace{6.5mm}
\raisebox{1.5mm}{$\frac{pq}{n-p}$}
\hspace{25mm}\\[3mm]
\hspace{79mm}
$\frac{1}{4n-4}$ \hspace{0mm}
$\frac{1}{2n-2}$ \hspace{1mm}
\dots
\hspace{1mm}
$\frac{1}{2n-2}$ \hspace{0mm}
$\frac{1}{4n-4}$ \hspace{17mm}\\
$\tilde{\pi}(x)$
\epsfig{figure=figure5.eps,width=46mm}
\hfill
$\tilde{\pi}(x) \tilde{K}_n(x,y)$ \hspace{1mm}
\epsfig{figure=figure5.eps,width=46mm}
\\[.5mm] 
\hspace{3mm}
$\frac{1}{2n-2}$ 
\hspace{2mm}$ \frac{1}{n-1}$
\hspace{2mm} $\dots$ \hspace{2mm} $ \frac{1}{n-1}$ 
\hspace{2mm}
$\frac{1}{2n-2}$
\hspace{27mm} \raisebox{1.5mm}{
$\frac{1}{4n-4}$} \hspace{6.5mm} \raisebox{1.5mm}{$\dots$}
\hspace{6.5mm}
\raisebox{1.5mm}{$\frac{1}{4n-4}$}
\hspace{15mm}\\[-1.5mm]
\caption{Stationary distributions and edge weights of $K$ and
$\tilde{K}$.
\label{symm3fig}}
\end{figure}

It is convenient to take 
\[
\tilde{K}_n  = 
\left[
\begin{array}{ccccc}
{1}/{2} & {1}/{2} & 0 &  \dots & 0 \\
1/4 & 1/2& 1/4&   \dots &0\\ 
\vdots  & \ddots &  \ddots & \ddots &  \vdots\\
 0 & \dots  & 1/4&  1/2 &1/4\\
 0 & \dots  & 0&  1/2 &1/2
\end{array}
\right],
\]
the ``lazy'' simple random walk on the $n$-path.
Substituting the stationary distributions and edge weights 
of $K_n$ and $\tilde{K}_n$ (as shown in Figure
\ref{symm3fig}) and  $\beta_1(\tilde{K}_n)=\frac{1}{2} +\frac{1}{2}
\cos \left(\frac{\pi }{n-1} \right)$ into~(\ref{eq:d-sc}) 
and~(\ref{eq:fill}) now gives
\begin{align*}
\beta_1(K_n) &= 1- \frac{pq\pi^2}{2(n-1)^2}
+\bigOfrac{1}{(n-1)^4}
\end{align*}
and
\[
\tvdist{P_n^{(n-1)t}(x_0, \cdot)}{\pi} 
\leq \frac{1}{2\sqrt{\pi(x_0)}}
\left( 1-\frac{pq \pi^2}{2(n-1)^2}
+\bigOfrac{1}{(n-1)^4}\right)^{t/2}.
\]
Since $\frac{q}{n-p} \leq 
\pi(x_0) \leq \frac{1}{n-p}$, we have shown that
$\bigO{n^2 \log n}$ steps suffice for $P_n^{n-1}$ to be within a fixed
distance of stationarity, while $\bigO{n^3 \log n}$
suffice for $P_n$. As can happen for bounds using only the second-largest
eigenvalue, these results are a factor of $\log n$ larger than necessary. 

\subsubsection*{Nash inequalities} 
Diaconis and Saloff-Coste \cite{diasal96b} bound convergence to
stationarity using Nash inequalities. Their results are often sharper than 
second-largest-eigenvalue bounds for slowly-converging chains. Here, we
extract only a few pieces of their analysis. 

Given an edge set $E$ on a state space
$S$, let
$d(x,y)$ be the length of a shortest path from $x$ to
$y$ via $E$. Let $B(x,r) = \{y: d(x,y)\leq
r\}$. Given a measure $\pi$ on $S$, let 
$V(x,r)=\sum_{y \in B(x,r)} \pi(y)$. Let $\gamma = 
\max_{x,y} d(x,y)$ be the diameter of $E$. 

Say
$(S,E,\pi)$ has $(A,d)$--{\em moderate growth} if
\[
V(x,r) \geq \frac{1}{A}
\left(\frac{r+1}{\gamma}\right)^d \quad
\text{for} \quad 0 \leq r \leq \gamma.\]
It is easily checked that the underlying graph and  stationary
distribution
 of $M(P_n^{n-1})$ have
$\left( \frac{1}{q}+\frac{1}{n-1}, 1\right)$--moderate
growth. 
\begin{Proposition}[Diaconis and Saloff-Coste
\cite{diasal96b}]
\label{pr:dsc}
Let $K$ be a Markov chain on a finite set $S$ with
stationary distribution $\pi$. Let $M(K)$ be the
multiplicative symmetrization of $K$, and $E$ the edge
set of the underlying graph of $M(K)$. Assume $(S,E,
\pi)$ has $(A,d)$--moderate growth. Let
$\{\gamma_{xy}\}$ be a collection of paths in $E$
joining each pair of states in $S$, and let
\[
a = \max_{
\begin{smallmatrix}
(x,y) \in E \\  0 \leq r \leq \gamma
\end{smallmatrix} }
\frac{2}{r^2\pi(x) M(K)(x,y)}
\left\{
\sum_{\begin{smallmatrix}
\gamma_{zw} \ni (x,y)\\ 
d(z,w) \leq r \end{smallmatrix}}
\abs{\gamma_{zw}} \frac{\pi(z)\pi(w)}{V(z,r)}
\right\}.
\]
Then 
\[
\tvdist{P^t(x_0, \cdot)}{\pi} \leq \frac{1}{2} 
a_1 e^{-m/(a \gamma^2)}, \quad \text{for} \quad
t = a\gamma^2 + m +1,
\]
when $m \geq 0$ and $a_1 = (e(1+d)A)^{1/2}
(4(2+d))^{d/4}$.
\end{Proposition}
Following example 5(a) of
\cite{diasal96b}, we take the paths $\{\gamma_{xy} \}$ in $K_n=M(P_n^{n-1})$ 
 to be geodesics and bound
$a$.    
From the values in Figure \ref{symm3fig}  we see that
\begin{align*}
\frac{2}{\pi(x)
M(K_n)(x,y)} \leq \frac{2(n-p)}{\min(q^2, pq)} \quad & \textrm{and} \quad
\pi(z)\pi(w)  \leq \frac{(1-pq)^2}{(n-p)^2}.
\end{align*}
The number of terms in the sum is at most 
$r(r+1)/2$ and $\abs{\gamma_{zw}} \leq r$. 
Finally, $\left(\frac{1}{q}+\frac{1}{n-1}, 1
\right)$--moderate growth implies that 
\begin{align*}
V(z,r) & \geq \frac{q(r+1)}{n-p}. 
\end{align*}
Combining these estimates yields
\begin{align*}
a & \leq \frac{(1-pq)^2}{q \min(q^2, pq)}.
\end{align*}
We can conclude that $\bigO{n^2}$ steps suffice for
the chain $P_n^{n-1}$ to be close to stationarity, and thus
$\bigO{n^3}$ suffice for the chain $P_n$.

\begin{Remark} Although it is difficult to compare the asymptotic statement of
Theorem~\ref{th:rate} and the direct inequalities of Proposition~\ref{pr:dsc},
the Nash inequality estimate of the lead term constant appears to worsen as $q$
decreases to 0. 

In order to force $\tvdist{P_n^{t(n-1)}(x_0,
\cdot)}{\pi}
\leq
\epsilon$ using Theorem~\ref{th:rate}, we need a
 $c$ such that 
\[
\left( \int_{\RZ} \abs{\thetafun{c}{x}-1} dx \right)^{1/2} < \epsilon.
\]
Then $t>\frac{c(n-p)^3}{pq(n-1)} = \frac{cn^2}{pq} +
\bigO{n}$ steps suffice.
In Proposition~\ref{pr:dsc}, we must have
\[
\quad m > \log \left(\frac{a_1}{\epsilon} \right)a
\gamma^2,
\]
which implies 
\[
t > \frac{(1-pq)^2}{q\min(q^2,pq)}\left(
1+\frac{1}{4} \log (48e^2) + \frac{1}{2} \log
\left(\frac{1}{q}\right) - \log \epsilon
\right) (n-1)^2 +1.
\]
As $q \rightarrow 0$, the
lead term is $\bigO{\frac{1}{q^3} \log \frac{1}{q}} n^2$, as opposed
to $\bigOfrac{1}{q} n^2$ for the asymptotic bound. 
\end{Remark}

\section*{Acknowledgements}
These results first appeared in the author's doctoral
dissertation~\cite{wilmer99}, which was supervised by Persi Diaconis. 
A suggestion of Ron
Graham's led to the examples studied here. The author is also grateful
to Jim Fill and Agoston Pisztora for helpful
conversations.

\end{document}